\documentclass[12pt]{article}
\usepackage{verbatim,amsmath,amsthm,amssymb,graphics,graphicx,color}
\newcommand{\finalmod}[1]{{\iftesting\color{red}\fi #1}}
\title{On the set-generic multiverse}

\renewcommand{\thefootnote}{}
\author{Sy David Friedman\,${}^\ast$, Saka\'e Fuchino\,${}^\dagger$ 
and Hiroshi Sakai\,${}^\ddagger$}

\newif\iftesting
\ifx\Testing\Kaba\else\testingtrue\fi
\newif\ifextended
\ifx\Extended\Kaba\else\extendedtrue\fi
\let\Label\label%
\iftesting
\def\label#1{\marginpar{{\renewcommand{\baselinestretch}{0.6}\tiny
			#1}}\Label{#1}\ignorespaces}%
\fi

\newcounter{frml}[section]
\def\thefrml{{\arabic{section}.\arabic{frml}}}
\def\frmlabel#1{\refstepcounter{frml}{\def\baka{#1}\ifx\baka\empty\else\Label{#1}\fi}%
{\rm({\thefrml})\hfill\hfill\hfill}}
\def\xitem[#1]{\item[\frmlabel{#1}]\mbox{}%
	\iftesting\marginpar{{\renewcommand{%
				\baselinestretch}{0.6}\tiny#1}}\fi\ignorespaces}
\def\xitemsub[#1]#2{\item[\frmlabel{#1}$_{#2}$]\mbox{}%
	\iftesting\marginpar{{\renewcommand{%
				\baselinestretch}{0.6}\tiny#1}}\fi\ignorespaces}
\def\xxitem[#1][#2]{\item[(\ref{#1}{\makebox[1.4ex][c]{#2}})]\mbox{}%
	\iftesting\marginpar{{\renewcommand{%
				\baselinestretch}{0.6}\tiny\{#1\}\{#2\}}}\fi\ignorespaces}
\def\itemof#1{{\bf(#1)}}
\def\xitemof#1{{\rm({\ref{#1}})}}

\newenvironment{xitemize}{\begin{list}{}{\parsep=0.5\smallskipamount%
			\itemindent=-0.0ex%
			\itemsep=0.5\smallskipamount\leftmargin=8ex\labelwidth=6ex\labelsep=1.4ex}}%
							 {\end{list}}
							 {\end{list}}

\def\assertof#1{{\rm (#1)}}
{\end{minipage}\end{trivlist}}
\newcommand{\bbd}[1]{{\mathbb{#1}}}
\makeatletter
\renewcommand{\section}{\@startsection {section}{1}{\z@ }{-3.5ex \@plus -1ex \@minus -.2ex}{2.3ex \@plus .2ex}{\normalfont \large \bfseries }}
\makeatother
\renewcommand{\baselinestretch}{1.2}

\newcommand{\ctentenc}{,\mbox{\hspace{0.04ex}{.}{.}{.}\hspace{0.1ex},\,}}

\newcommand{\xmbox}[1]{ $\relax{\rm #1}\relax$ }
\newcommand{\llor}{{\bigvee\hspace{-1.5ex}\bigvee}\rule[-0.8ex]{0cm}{1ex}}
\newcommand{\lland}{{\bigwedge\hspace{-1.5ex}\bigwedge}\rule[-0.8ex]{0cm}{1ex}}
\newcommand{\cardof}[1]{\mathopen{|\,}#1\mathclose{\,|}}

\newcommand{\setof}[2]{\{#1\,:\,#2\}}
\newcommand{\ssetof}[1]{\{#1\}}

\newcommand{\seqof}[2]{\langle#1\,:\,#2\rangle}
\newcommand{\pairof}[1]{\langle#1\rangle}

\newcommand{\psof}[1]{{\mathcal P}\/(#1)}
\newcommand{\Pof}[2]{{{\mathcal P}_{\!#1}#2}}

\newcommand{\mapping}[3]{#1:#2\rightarrow #3}

\newcommand{\fnsp}[2]{\mbox{}^{#1^{\mbox{}\!}}#2}
\newcommand{\circleq}{\mathrel{{\leq}%
		\hspace{-0.86ex}{\lower-0.5ex\hbox{$\scriptscriptstyle\circ$}}}}

\newcommand{\forces}[2]{\,\|\hspace{-.35ex}\mbox{\sf--}_{\,#1\,}%
\mbox{\rm``}\,#2\,\mbox{\rm''}}

\newcommand{\utildeT}[1]{%
	\hspace{0.25ex}\hbox to 0pt%
				 {$\mathop{#1}\limits_{\raise0.24ex\hbox{$\scriptstyle\sim$}}$\hss}%
		\relax\phantom{\underline{#1\hspace{0.2ex}}}}
\newcommand{\utildeS}[1]{%
	\hspace{0.25ex}\hbox to 0pt%
				 {$\mathop{\scriptstyle #1}\limits_{%
						 \raise0.44ex\hbox{$\scriptscriptstyle\sim$}}$\hss}%
		\hspace{0.25ex}\relax\phantom{\underline{#1\hspace{0.25ex}}}}
\newcommand{\utildeSS}[1]{%
	\hbox to 0pt{$\mathop{\scriptscriptstyle #1}%
		\limits_{\raise0.44ex\hbox{$\scriptscriptstyle\sim$}}$\hss}%
		\relax\phantom{\underline{#1\hspace{0.25ex}}}}

\newcommand{\On}{{\rm On}}

\newcommand{\dom}{\mathop{\rm dom}}

\newcommand{\range}{\mathop{\rm rng}}

\newcommand{\Fn}{\mathop{\rm Fn}}

\newcommand{\GCH}{{\rm GCH}}

\newcommand{\ZF}{{\rm ZF}}
\newcommand{\AC}{{\rm AC}}

\newcommand{\ZFC}{{\rm ZFC}}
\newcommand{\NBG}{{\rm NBG}}
\newcommand{\MZF}{{\rm MZF}}
\newcommand{\MZFC}{{\rm MZFC}}
\newcommand{\LZF}{{\mathcal L}_{\rm ZF}}

\newtheorem{Thm}{Theorem}[section]

\newtheorem{Cor}[Thm]{Corollary}

\newtheorem{Prop}[Thm]{Proposition}
\newtheorem{Lemma}[Thm]{Lemma}

\newtheorem{Claim}{Claim}[Thm]

\newcommand{\Thmof}[1]{{\rm Theorem \ref{#1}}}
\newcommand{\bfThmof}[1]{{\bf Theorem \ref{#1}}}

\newcommand{\Corof}[1]{{\rm Corollary \ref{#1}}}

\newcommand{\Lemmaof}[1]{{\rm Lemma \ref{#1}}}
\newcommand{\bfLemmaof}[1]{{\bf Lemma \ref{#1}}}

\newcommand{\Claimof}[1]{{\rm Claim \rm\ref{#1}}}
\newcommand{\bfClaimof}[1]{{\bf Claim \ref{#1}}}

\newcommand{\sectionof}[1]{{\rm Section \ref{#1}}}

\newcommand{\Corabove}{{\rm Corollary \number\theThm}}

\newcommand{\prf}{\noindent{\bf Proof.\ }\ignorespaces}
\newcommand{\prfof}[1]{\noindent{\bf Proof of #1:\,\ }\ignorespaces}
\newcommand{\prfofClaim}{\noindent\raisebox{-.4ex}{\Large $\vdash$\ \ }}

\newsavebox{\qedbox}\sbox{\qedbox}{
{\unitlength=0.065mm \begin{picture}(40,60)
\put(0,0){\framebox(30,44)[lb]{}}
\put(30,-6){\rule{6\unitlength}{44\unitlength}}
\put(10,-6.2){\rule{26\unitlength}{6\unitlength}}
\end{picture}}}
\renewcommand{\qed}{\mbox{}\hfill\usebox{\qedbox}}
\newcommand{\smallqed}%
{\mbox{}\smallskip\hfill\raisebox{-.4ex}{\Large $\dashv$}\\}
\newcommand{\qedof}[1]%
{\mbox{} \hspace*{\fill}{\usebox{\qedbox}{\rm~(#1)}}%
\mbox{}}
\newcommand{\Qedof}[1]%
{\mbox{} \hspace*{\fill}{\usebox{\qedbox}%
{\rm~(#1~\number\theThm)}}}
\newcommand{\QedAof}[1]%
{\mbox{} \hspace*{\fill}{\usebox{\qedbox}%
{\rm~(#1~A.\number\theThmA)}}}
\newcommand{\qedofThm}{\Qedof{Theorem}}

\newcommand{\qedofProp}{\Qedof{Proposition}}
\newcommand{\qedofLemma}{\Qedof{Lemma}}

\newcommand{\qedskip}{\medskip\smallskip}
\newcommand{\qedofClaim}%
{\mbox{}\hfill\raisebox{-.4ex}{\Large $\dashv$ }\nolinebreak%
\mbox{\rm~({\rm Claim}~\rm\number\theClaim)}}
\newcommand{\qedofClaimA}%
{\mbox{}\hfill\raisebox{-.4ex}{\Large $\dashv$ }\nolinebreak%
\mbox{\rm~({\rm Claim}~\rm{A.}\number\theClaimA)}}
\newcommand{\qedofSubclaim}%
{\mbox{}\hfill\raisebox{-.4ex}{\Large $\dashv$ }\nolinebreak%
\mbox{\rm~({\rm Subclaim}~\rm\number\theSubclaim)}}
\newcommand{\calA}{{\mathcal A}}

\newcommand{\calH}{{\mathcal H}}

\newcommand{\calL}{{\mathcal L}}
\newcommand{\calM}{{\mathcal M}}

\newcommand{\calP}{{\mathcal P}}

\newcommand{\dota}{\dot{a}}

\newcommand{\dotf}{\dot{f}}

\newcommand{\dotx}{\dot{x}}
\newcommand{\dotA}{\dot{A}}

\newcommand{\dotG}{\dot{G}}

\newcommand{\dotS}{\dot{S}}
\newcommand{\dotT}{\dot{T}}
\newcommand{\bbB}{\bbd{B}}
\newcommand{\sfS}{{\sf S}}
\newcommand{\sfT}{{\sf T}}
\newcommand{\po}{{partial ordering}}

\newcommand{\poP}{\bbd{P}}

\newcommand{\bbbone}{{\mathchoice {\rm 1\mskip-4mu l} {\rm 1\mskip-4mu l}
{\rm 1\mskip-4.5mu l} {\rm 1\mskip-5mu l}}}

\newcommand{\BAB}{\bbd{B}}
\newcommand{\st}{such that}

\newcommand{\Wolog}{Without loss of generality}
\newcommand{\wrt}{with respect to}

\newcount\minute	
\newcount\hour		
\newcount\hourMins  
\def\now%
{
  \minute=\time    
  \hour=\time \divide \hour by 60 
  \hourMins=\hour \multiply\hourMins by 60
  \advance\minute by -\hourMins 
  \zeroPadTwo{\the\hour}:\zeroPadTwo{\the\minute}%
}
\def\zeroPadTwo#1%
{
  \ifnum #1<10 0\fi    
  #1
}

\date{}
\begin{document}
\footnotetext{{\it Date:} \today\ (\now\ifx\kanjiskip\nanjamonja{}\else\ JST\fi)
\vspace{-1\smallskipamount}}
\footnotetext{
\vspace{-1\smallskipamount}}
\footnotetext{
  \vspace{-1\smallskipamount}} 
\renewcommand{\thefootnote}{$\ast$}
\footnotetext{
  Kurt G\"odel Research Center for Mathematical Logic,
  University of Vienna, Vienna, Austria, E-Mail: {\tt sdf@logic.univie.ac.at}\quad
  The first author would like to thank the FWF (Austrian Science Fund)
  for its support through project number P 28420.
  
  $\ \!\!{}^\dagger$Graduate School of System Informatics, Kobe University, 
  Kobe, Japan, E-Mail: {\tt fuchino@diamond.kobe-u.ac.jp}\quad
  The second author is supported by Grant-in-Aid for
  Scientific Research (C) No.\ 21540150 and Grant-in-Aid for Exploratory 
  Research No.\ 26610040 of the Ministry of Education,
  Culture, Sports, Science and Technology Japan (MEXT).

  $\ \!\!{}^\ddagger$Graduate School of System Informatics, Kobe University,  
  Kobe, Japan, E-Mail: {\tt hsakai@people.kobe-u.ac.jp}\quad
  The third  author is supported by Grant-in-Aid for Young Scientists (B) No.\ 23740076
  of the Ministry of Education, Culture, Sports, Science and Technology Japan (MEXT). 
}
\renewcommand{\thefootnote}{}
\footnotetext{The second author would like to thank Toshimichi Usuba for some 
  valuable comments. \finalmod{The authors also would like to thank the anonymous 
    referee for many valuable comments and suggestions.}}
\maketitle


\begin{abstract}
The forcing method is a powerful tool to prove the consistency of set-theoretic 
assertions relative to the consistency of the axioms of set theory. Laver's 
theorem and Bukovsk\'y's theorem assert that set-generic extensions of a given 
ground model constitute a quite reasonable and  
sufficiently general class of standard models of set-theory. 

In sections 2 and 3 of this note, we give a proof of
Bukovsky's theorem in a modern setting (for another proof
of this theorem see \cite{bukovsky2}). In section 4 we check that the
multiverse of set-generic extensions can be treated as a
collection of countable transitive models in a conservative
extension of \ZFC. The last section then deals with the
problem of the existence of infinitely-many independent
buttons, which arose in the modal-theoretic approach to the
set-generic multiverse by J.\,Hamkins and B.\,Loewe \cite{hamkins-loewe}.\end{abstract}
\renewcommand{\thefootnote}{(\arabic{footnote})\,}
\section{The category of forcing extensions as the set-theoretic multiverse}

The forcing method is a powerful tool to prove the consistency of set-theoretic 
(i.e., mathematical) assertions relative to (the consistency of) the axioms of 
set theory. If a sentence  
$\sigma$ in the language $\LZF$ of set theory is proved to be relatively 
consistent with the axioms of set theory (\ZFC) by some 
forcing argument then it is so in the sense of the strictly finitist standpoint of Hilbert: 
the forcing proof can be recast into an algorithm $\calA$ \st, if a 
formal proof $\calP$ of a contradiction from \ZFC\ $+$ $\sigma$ is ever given, then we 
can transform $\calP$ with the help of $\calA$ to another proof of a 
contradiction from \ZFC\ or even \ZF\ alone. 

The ``working set-theorists'' however prefer to see their forcing arguments not as 
mere discussions concerning manipulations of formulas in a formal system but rather concerning 
the ``real'' 
mathematical universe in which they ``live''. Forcing for them is thus a method of 
extending the universe 
of set theory where they originally ``live'' (the ground model, usually denoted as
``$V$'') to many (actually more than class many in the sense of $V$) different 
models of set theory  
called generic extensions of $V$. Actually, a family of generic extensions is constructed 
for certain $V$-definable partial orderings $\poP$. Each such generic extension is 
obtained first by fixing a so-called generic filter $G$ which is 
a filter over $\poP$, sitting outside $V$ with a ``generic'' sort of 
transcendence over $V$, and then by adding $G$ to $V$ to generate a new structure 
--- the generic extension $V[G]$ of $V$ --- which is also a model of \ZFC. Often this 
process of taking generic extensions over some model of set theory 
is even repeated transfinitely-many times.  As a result, a set-theorist performing forcing 
constructions is seen to live in   
many different models of set theory simultaneously. This is manifested in many 
technical expositions of forcing where the reader very often finds narratives  
beginning with phrases like: ``Working in $V[G]$, \ldots'', ``Let $\alpha<\kappa$ 
be \st\ $x$ is in the $\alpha$ th intermediate model $V[G_\alpha]$ and \ldots'', 
``Now returning to $V$, \ldots'', etc., etc.

Although this ``multiverse'' view
of forcing is in a sense merely a modus loquendi, it is worthwhile to study the possible 
pictures of this multiverse per se. Some initial moves in this 
direction have been taken e.g. in 
\cite{tatiana-sy-1}, \cite{tatiana-sy-2}, \cite{sy-strict-gen}, \cite{sy-book}, \cite{sy-imh},
\cite{fuhare}, \cite{hamkins-multiverse}, \cite{hamkins-loewe}, \cite{usuba}, 
\cite{woodin-realm} etc.  
The term ``multiverse'' probably originated in work of Woodin in which he considered the 
``set-generic multiverse'', the ``class'' of set-theoretic universes which forms 
the closure of the 
given initial universe $V$ under set-generic extension and set-generic 
ground models. Sometimes we also have to consider the constellations of 
the set-generic multiverse where $V$ cannot be reconstructed as a set-generic 
extension of some of or even any of the proper inner models of $V$. To deal with such 
cases it is more convenient to consider the expanded generic multiverse where we 
also assume that the multiverse is also closed under the construction of 
definable inner models. 

The set-generic universe should be distinguished from the ``class-generic 
multiverse'', defined in the 
same way but with respect to class-forcing extensions and ground models, as well as inner models
of class-generic extensions that are not themselves class-generic (see \cite{sy-strict-gen}). It is 
even possible to go beyond class-forcing by considering forcings whose conditions are classes, so-called
hyperclass forcings (see \cite{sy-book}). The broadest point of view with regard to the multiverse
is expressed in \cite{sy-imh}, where the ``hyperuniverse'' is taken to consist of \emph{all} universes
which share the same ordinals as the initial universe (which is taken to be 
countable to facilitate 
the construction of new universes). The hyperuniverse is closed under all notions of forcing.

In this article we restrict our attention to the set-generic multiverse.
The well-posedness of questions regarding the set-generic multiverse is
established by the theorems of Laver and Bukovsk\'y which we 
discuss in Section \ref{laver-bukovsky}. These theorems show that the set-generic 
extensions and set-generic ground models of a given universe represent a ``class'' of
models with a natural characterization.

The straightforward formulation of the set-generic multiverse requires the notion of ``class'' 
of classes which cannot be treated in the usual framework of \ZF\ set theory,
but, as emphasized at the beginning, theorems about the set-generic multiverse are actually 
meta-theorems about \ZFC. However we can also consider a theory which is a 
conservative extension of \ZFC\ in which set-generic extensions and set-generic 
ground models are real objects in the theory and the set-generic multiverse a definable class. 
In Section \ref{multi-sys}, we 
consider such a system and show that it is a conservative extension of 
\ZFC. 

The multiverse view sometimes highlights problems which would never have been
asked in the conventional context of forcing constructions (see \cite{hamkins-multiverse}). 
As one such example we consider in Section \ref{ind-buttons} the problem of the existence of 
infinitely many independent buttons (in the sense of \cite{hamkins-loewe}). 

\section{Laver's theorem and Bukovsk\'y's theorem}
\label{laver-bukovsky}
In the forcing language, we often have to express that a certain set is already 
in the ground model, e.g.\ in a statement like: 
$p\forces{\poP}{\ldots\ \dotx\mbox{ is in }V\mbox{ and }\ldots}$. In such 
situations we can always find a large enough ordinal $\xi$ \st\ the set in question 
should be found in that level of the cumulative hierarchy in the ground model.  So we can 
reformulate a statement like the one above into something like 
$p\forces{\poP}{\ldots\ \dotx\in\check{V}_\xi\mbox{ and }\ldots}$ which is a 
legitimate expression in the forcing language.

This might be one of the reasons why it is proved only quite recently that the 
ground model is always definable in an arbitrary set-generic extension:
\begin{Thm}[R.\,Laver, \cite{laver}, H.\,Woodin \cite{woodin}]
  \label{laver-thm}
  There is a formula $\varphi^*(x,y)$ in 
  $\calL_\ZF$ \st, for any transitive model $V$ of \ZFC\ and set-generic extension 
  $V[G]$ of $V$ there is $a\in V$ \st, for any $b\in V[G]$ 
  \[b\in V\ \ \Leftrightarrow\ \ V[G]\models\varphi^*(a,b).  \vspace{-6ex}
  \]\noindent
  \mbox{}\qed
\end{Thm}

An important corollary of Laver's theorem is that a countable transitive model of $\ZFC$ 
can have at most countably many ground models for set forcing.

%
%

Bukovsk\'y's theorem gives a natural characterization of inner models $M$ of $V$ 
\st\ $V$ is a set-generic extension of $M$\footnote{In the terminology of 
  \cite{fuhare}, $M$ is a ground of $V$.}. Note that, by Laver's theorem 
\Thmof{laver-thm}, such an $M$ is then definable in $V$.  However the inner model
$M$ of $V$ may be introduced as a class in the sense of von 
Neumann-Bernays-G\"odel class theory (NBG) and in such a situation the 
definability of $M$ 
in $V$ may not be immediately clear. 

Let us begin with the following observation concerning $\kappa$-c.c.\ generic extensions. 
We shall call a partial ordering \emph{atomless} if each element of it has at 
least two extensions which are incompatible with each other.

\begin{Lemma}
  \label{kappa-c-c}
  Let $\kappa$ be a regular uncountable cardinal. If\/ $\poP$ is a $\kappa$-c.c.\ 
  atomless \po, then $\poP$ adds a new subset of $2^{<\kappa}$. 
\end{Lemma}
\prf \Wolog, we may assume that $\poP$ consists of the positive elements of a
$\kappa$-c.c.\  
atomless complete Boolean algebra. Note that $\poP$ adds new subset of $\On$ 
since $\poP$ adds a new set (e.g.\ the 
$(V,\poP)$-generic set).  
Suppose that $\dotS$ is a $\poP$-name of a new subset 
of $\On$. Let $\theta$ be a sufficiently large regular cardinal and let 
$M\prec\calH(\theta)$ be \st\
\begin{xitemize}
\xitem[buk-0] $\cardof{M}\leq 2^{<\kappa}$;
\xitem[buk-1] $\fnsp{<\kappa}{M}\subseteq M$ and 
\xitem[buk-2] $\poP$, $\dotS$, $\kappa\in M$.
\end{xitemize}
Let $\dotT$ be a $\poP$-name \st $\forces{\poP}{\dotT=\dotS\cap M}$. By 
\xitemof{buk-0}, it is enough 
to show the following, where $V$ denotes the ground model:
\begin{Claim}
  $\forces{\poP}{\dotT\not\in V}$. 
\end{Claim}
\prfofClaim
Otherwise there would be $p\in\poP$ and $T\in V$, $T\subseteq\On$ \st\
\begin{xitemize}
\xitem[] 
  $p\forces{\poP}{\dotT=\check{T}}$. 
\end{xitemize}
We show in the following that then we can construct a strictly decreasing 
sequence $\seqof{q_\alpha}{\alpha<\kappa}$ in $\poP\cap M$ \st\
\begin{xitemize}
\xitem[buk-3] 
  $p\leq_\poP q_\alpha$ for all $\alpha<\kappa$. 
\end{xitemize}
But since
$\setof{q_{\alpha}\cdot-q_{\alpha+1}}{\alpha<\kappa}$ is then a pairwise disjoint 
subset of $\poP$, this contradicts the $\kappa$-c.c. of $\poP$. 

Suppose that $\seqof{q_\alpha}{\alpha<\delta}$ for some $\delta<\kappa$ has been 
constructed. If $\delta$ is a limit, let $q_\delta=\prod_{\alpha<\delta}q_\alpha$. 
Then we have $p\leq_\poP q_\delta$ and $q_\delta\leq_\poP q_\alpha$ for all $\alpha<\delta$.
Since
$\seqof{q_\alpha}{\alpha<\delta}\in M$ by \xitemof{buk-1}, we also have
$q_\delta\in M$. 

If $\delta=\beta+1$, then, since
$M\models\xmbox{``}q_\beta\xmbox{ does not decide }\dotS\xmbox{''}$ by the 
elementarity of $M$, there are $\xi\in\On\cap M$ and  $q$,
$q'\in \poP\cap M$ with $q$, $q'\leq_\poP q_\beta$ 
\st\ $q\forces{\poP}{\xi\in\dotS}$ 
and $q'\forces{\poP}{\xi\not\in\dotS}$. At least one of them, say $q$, must be 
incompatible with $p$. Then $q_\delta=q_\beta\cdot-q$ is as desired.
\qedofClaim\\
\qedofLemma\qedskip

Note that, translated into the language of complete Boolean algebras, the lemma 
above just asserts that no $\kappa$-c.c.\ atomless Boolean algebra $\bbB$ is 
$(2^{<\kappa},2)$-distributive. 

Suppose now that we work in NBG, $V$ is a transitive model of \ZFC\ and $M$ an 
inner model of \ZFC\  
in $V$ (that is $M$ is a transitive class $\subseteq V$ 
with $(M,\in)\models\ZFC$). For a regular uncountable cardinal $\kappa$ in $M$, 
we say  
that $M$ {\it $\kappa$-globally covers $V$\/} if for every function $f$ (in $V$) with
$\dom(f)\in M$ and $\range(f)\subseteq M$, there is a function $g\in M$ with
$\dom(g)=\dom(f)$ \st\ 
$f(i)\in g(i)$ and $M\models\cardof{g(i)}<\kappa$ for all
$i\in\dom(f)$. 

\begin{Thm}[L.\ Bukovsk\'y, \cite{bukovsky}\footnotemark]
  \label{bukovsky-thm}
  Suppose that $V$ is a transitive model of\/ \ZFC, $M\subseteq V$ an inner model of\/ 
  \ZFC\ and $\kappa$ is a regular uncountable cardinal in $M$
  . 
  Then $M$ $\kappa$-globally covers $V$ if and only if\/ $V$ is  
  a $\kappa$-c.c.\  set-generic extension of $M$.
\end{Thm}
\footnotetext{Tadatoshi Miyamoto told us that 
  James Baumgartner independently proved this theorem in an unpublished note using 
  infinitary logic.}%

\finalmod{As the referee of the paper points out, this theorem can be formulated 
  more naturally in the von Neumann-Bernays-G\"odel class theory (NBG) 
  since in the framework of \ZFC\ this theorem can only be formulated as a 
  meta-theorem, that is, as a collection of theorems consisting corresponding 
  statements for each formula which might define an inner model $M$. }\qedskip

\prfof{\bfThmof{bukovsky-thm}} If $V$ is a $\kappa$-c.c.\ set-generic extension 
of $M$, say by a  
\po\ $\poP\in M$ with $M\models$``$\poP$ has the $\kappa$-c.c.'', then it is clear that
$M$ $\kappa$-globally covers $V$ (for $f$ as above, let 
$\dotf\in M$ be a $\poP$-name of $f$ and $g$ be defined by letting $g(\alpha)$ to 
be the set of all possible values $\dotf(\alpha)$ may take).  

The proof of the converse is done via the following \Lemmaof{buk-main-cl}. 
Note that, by Grigorieff's theorem (see \Corof{grigorieff} below), the statement 
of this Lemma is a consequence of Bukovsk\'y's theorem: 
\begin{Lemma}
  \label{buk-main-cl}
  Suppose that $M$ is an inner model of a transitive model $V$ of \ZFC\ \st\ $M$
  $\kappa$-globally covers $V$ for some $\kappa$ regular uncountable in $M$.
  Then for any $A\in V$, $A\subseteq\On$, $M[A]$ is\footnote{%
    \Label{M[A]}
    $M[A]$ may be defined by
    $M[A]=\bigcup_{\alpha\in\On}L(V_\alpha^M\cup\ssetof{A})$. 
    \finalmod{$M[A]$ is a model of \ZF: this can be seen easily e.g.\ by applying 
  Theorem 13.9 in \cite{millennium-book}. If $M$ also satisfies \AC\ then $M[A]$ 
  satisfies \AC\ as well since, in this case, it is easy to see that a
  well-ordering of $(V_\alpha)^M\cup\ssetof{A}$ belongs to $M[A]$ for all
  $\alpha\in\On$. }} a $\kappa$-c.c.\ set-generic extension  
  of $M$. 
\end{Lemma}

\finalmod{Note that it can happen easily that $M[A]$ is not a set 
  generic extrension of $M$. For example, $0^{\#}$ exists and $M=L$, then
  $M[0^{\#}]$ is not a set-generic extension of $M$.
}

We first show that \Thmof{bukovsky-thm} follows from \Lemmaof{buk-main-cl}. 
Assume that $M$
$\kappa$-globally covers $V$.  We have to show that $V$ is a $\kappa$-c.c.\ set-generic 
extension of $M$.
In $V$, let $\lambda$ be a regular cardinal \st\ $\lambda^{<\kappa}=\lambda$ and 
$A\subseteq\On$ be a set \st\
\begin{xitemize}
\xitem[buk-4] 
  $(\psof{\lambda})^{M[A]}=(\psof{\lambda})^V$. 
\end{xitemize}
Then, by 
\Lemmaof{buk-main-cl}, $M[A]$ is a $\kappa$-c.c.\ generic extension of $M$ and 
hence we have $M[A]\models\mbox{``}\kappa\xmbox{ is a regular cardinal''}$. 
Actually we have $M[A]=V$. Otherwise there would be a $B\in V\setminus M[A]$ 
with $B\subseteq\On$. Since $M[A]$ $\kappa$-globally covers $M[A][B]$, we may 
apply \Lemmaof{buk-main-cl} on this pair and conclude that $M[A][B]$ is a
(non trivial) $\kappa$-c.c.\ generic extension of $M[A]$. By \Lemmaof{kappa-c-c}, 
there is a new element of $\psof{(2^{<\kappa})^{M[A]}}\subseteq\psof{\lambda}$ in
$M[A][B]$.
 But this is a contradiction to \xitemof{buk-4}.\qedof{\Thmof{bukovsky-thm}}\qedskip

\prfof{\bfLemmaof{buk-main-cl}}
We work in $M$ and construct a $\kappa$-c.c.\ \po\ $\poP$ \st\ $M[A]$ 
is a $\poP$-generic extension over $M$. 

\finalmod{%
  Let $\mu\in\On$ be \st\ $A\subseteq\mu$ and  
  let 
  $\calL_{\infty}(\mu)$ be the infinitary sentential logic with atomic sentences
  \begin{xitemize}
  \xitem[buk-5] 
    ``$\alpha\in\dotA$'' for $\alpha\in\mu$ 
  \end{xitemize}
  and the class of sentences closed under $\neg$ and $\llor$ 
  where $\neg$ is to be applied to a formula and $\llor$ to an arbitrary set of 
  formulas. }
\finalmod{%
  To be specific let us assume that the atomic sentences  ``$\alpha\in\dotA$'' for
  $\alpha\in\mu$ are coded by the sets $\pairof{\alpha,0}$ for $\alpha\in\mu$, the negation 
  $\neg\varphi$ by $\pairof{\varphi,1}$ and the infinitary 
  disjunction $\llor\Phi$ by $\pairof{\Phi,2}$. }
We regard the usual disjunction $\lor$ of two formulas as a special case of 
$\llor$ and other logical connectives like ``$\lland$'', ``$\land$'', ``$\rightarrow$'' as  
being introduced as abbreviations of usual combinations of $\neg$ and
$\llor$. 
For a sentence $\varphi\in\calL_{\infty}(\mu)$ and $B\subseteq\mu$, we write 
$B\models\varphi$ when $\varphi$ holds if each atomic sentence of the form
``$\alpha\in\dotA$'' in $\varphi$ is interpreted by ``$\alpha\in B$'' and logical 
connectives in $\varphi$ are interpreted in canonical way. For a set $\Gamma$ of 
sentences, we write $B\models\Gamma$ if $B\models\psi$ for all $\psi\in\Gamma$. 
For $\Gamma\subseteq\calL_{\infty}(\mu)$ and $\varphi$, we write 
$\Gamma\models\varphi$ if $B\models\Gamma$ implies $B\models\varphi$ for all
$B\subseteq\mu$ (in $V$). 


Let $\vdash$ be a notion of provability for $\calL_{\infty}(\mu)$ in some 
logical system which is correct (i.e.\ $\Gamma\vdash\varphi$ always implies
$\Gamma\models\varphi$)\footnote{More precisely, we assume that 
  \ZFC\ proves the correctness of $\vdash$.}, upward absolute (i.e.\
$M\subseteq N$ and $M\models\mbox{``}\Gamma\vdash\varphi\mbox{''}$ always imply 
$N\models\mbox{``}\Gamma\vdash\varphi\mbox{''}$ for any transitive 
models $M$, $N$ of \ZF) and  
sufficiently strong (so that all the arguments used below work for this $\vdash$). 
In \sectionof{ded-sys} we introduce one such deductive system (as well as an alternative
approach without using such a deduction system, based on L\'evy Absoluteness).

\finalmod{Let $\lambda=\max\ssetof{\kappa,\mu^+}$ and 
  $\calL_{\lambda}(\mu)=\calL_\infty(\mu)\cap (V_\lambda)^M$.}
Let $f\in V$ be a mapping 
$\mapping{f}{\big(\psof{\calL_{\lambda}(\mu)}\big)^M\setminus\ssetof{\emptyset}}{%
  \big(\calL_{\lambda}(\mu)\big)^M}$ 
\st, for any
$\Gamma\in \big(\psof{\calL_{\lambda}(\mu)}\big)^M\setminus\ssetof{\emptyset}$, we have
$f(\Gamma)\in\Gamma$ and $A\models f(\Gamma)$ if $A\models\llor\Gamma$. 
Since $M$ $\kappa$-globally covers $V$, there is a $g\in M$ with 
$\mapping{g}{\big(\psof{\calL_{\lambda}(\mu)}\big)^M\setminus\ssetof{\emptyset}}{%
  \Pof{<\kappa}{\big(\calL_{\lambda}(\mu)\big)}^M}$ \st\
$f(\Gamma)\in g(\Gamma)\subseteq\Gamma$ for all
$\Gamma\in (\psof{\calL_{\lambda}(\mu)})^M\setminus\ssetof{\emptyset}$. 

In $M$, let 
\begin{xitemize}
\xitem[] $T=\setof{\llor\Gamma\rightarrow\llor g(\Gamma)}{\Gamma\in
  \psof{\calL_{\lambda}(\mu)}\setminus\ssetof{\emptyset}}$. 
\end{xitemize}
Note that $M[A]\models\mbox{``}A\models T\mbox{''}$. It follows that $T$ is 
consistent \wrt\ our deduction system (in $V$).
In $M$, let 
\begin{xitemize}
\xitem[buk-6] 
  $\poP=\setof{\varphi\in\calL_{\lambda}(\mu)}{T\not\vdash\neg\varphi}$
\end{xitemize}
and for $\varphi$, $\psi\in\poP$, let
\begin{xitemize}
\xitem[buk-7] 
  $\varphi\leq_\poP\psi$\ \ $\Leftrightarrow$\ \ $T\vdash\varphi\rightarrow\psi$. 
\end{xitemize}

\begin{Claim}
  \label{buk-poP}
  For $\varphi\in\calL_{\lambda}(\mu)$, if $A\models\varphi$ then we have 
  $\varphi\in\poP$. 
  In particular, 
  $\mbox{``}\alpha\in\dotA\mbox{''}\in\poP$ for all $\alpha\in A$ and 
  $\mbox{``}\neg(\alpha\in\dotA)\mbox{''}\in\poP$ for all $\alpha\in\mu\setminus A$. 
\end{Claim}
\prfofClaim
Suppose $A\models\varphi$. We have to show $T\not\vdash\neg\varphi$: If 
$T\vdash\neg\varphi$ in $M$, 
then we would have
$V\models\mbox{``}T\vdash\neg\varphi\mbox{''}$. Since
$A\models T$ in $V$, it follows 
that $A\models\neg\varphi$.  This is a contradiction. 
\qedofClaim\qedskip

\begin{Claim}
  \label{buk-comp}
  For $\varphi$, $\psi\in\poP$, $\varphi$ and $\psi$ are compatible if and only if 
  \begin{xitemize}
  \xitem[buk-8] 
    $T\not\vdash\neg(\varphi\land\psi)$.  
  \end{xitemize}
\end{Claim}
Note that \xitemof{buk-8} is equivalent to 
\begin{xitemize}
\xitem[buk-8-0] 
  $T\not\vdash\neg\varphi\lor\neg\psi$
    \ \ \ (\,$\Leftrightarrow$\ \ \  $T\not\vdash\varphi\rightarrow\neg\psi$). 
\end{xitemize}
\prfofClaim
Suppose that $\varphi$, $\psi\in\poP$ are compatible. By the definition of
$\leq_\poP$ this means that there is  
$\eta\in\poP$ \st\ $T\vdash\eta\rightarrow\varphi$ 
and $T\vdash\eta\rightarrow\psi$. For this $\eta$ we have
$T\vdash\eta\rightarrow(\varphi\land\psi)$. Since $T\not\vdash\neg\eta$ by the 
consistency of $T$, it 
follows that $T\not\vdash\neg(\varphi\land\psi)$. 

Conversely if $T\not\vdash\neg(\varphi\land\psi)$. Then 
$(\varphi\land\psi)\in\poP$. Since $T\vdash(\varphi\land\psi)\rightarrow\varphi$ 
and $T\vdash(\varphi\land\psi)\rightarrow\psi$, we have 
$(\varphi\land\psi)\leq_\poP\varphi$ and $(\varphi\land\psi)\leq_\poP\psi$. Thus 
$\varphi$ and $\psi$ are compatible \wrt\ $\leq_\poP$. 
\qedofClaim

\begin{Claim}
  \label{buk-kappa-cc}
  $\poP$ has the $\kappa$-c.c.
\end{Claim}
\prfofClaim
Suppose that $\Gamma\subseteq\poP$ is an antichain. Since
$\cardof{g(\Gamma)}<\kappa$, it is enough 
to show that $g(\Gamma)=\Gamma$. Suppose otherwise and let
$\varphi_0\in\Gamma\setminus g(\Gamma)$. Since
$\mbox{``}\llor\Gamma\rightarrow\llor g(\Gamma)\mbox{''}\in T$ and
$\vdash\varphi_0\rightarrow\llor\Gamma$, we have
\begin{xitemize}
\xitem[buk-9] 
  $T\vdash\varphi_0\rightarrow\llor g(\Gamma)$.
\end{xitemize}
It follows that there is $\varphi\in g(\Gamma)$ \st\ $\varphi_0$ and $\varphi$ 
are compatible. This is because otherwise we would have
$T\vdash\varphi_0\rightarrow\neg\varphi$ for all $\varphi\in g(\Gamma)$ by 
\Claimof{buk-comp}. Hence
$T\vdash\varphi_0\rightarrow\lland\setof{\neg\varphi}{\varphi\in g(\Gamma)}$ which is 
equivalent to $T\vdash\varphi_0\rightarrow\neg\llor g(\Gamma)$. From this and 
\xitemof{buk-9}, it follows that 
$T\vdash\neg\varphi_0$. But this is a contradiction to the 
assumption that $\varphi_0\in\poP$. 

Now, since $\Gamma$ is pairwise incompatible, it follows that
$\varphi_0=\varphi\in g(\Gamma)$. This is a contradiction to the choice of
$\varphi_0$. \qedofClaim\qedskip 

In $V$, let $G(A)=\setof{\varphi\in\poP}{A\models\varphi}$. By \Claimof{buk-poP}, 
we have $G(A)=\setof{\varphi\in\calL_{\lambda}(\mu)}{A\models\varphi}$ and 
$A$ 
is definable from $G(A)$ over $M$ as
$\setof{\alpha\in\mu}{\mbox{``}\alpha\in\dot{A}\mbox{''}\in G(A)}$. Thus we have
  $M[G(A)]=M[A]$.  

Hence the following two Claims prove our Lemma: 
\begin{Claim}
  $G(A)$ is a filter in $\poP$.
\end{Claim}
\prfofClaim
Suppose that $\varphi\in G(A)$ and $\varphi\leq_\poP\psi$. Since this means that 
$A\models\varphi$ and $T\vdash\varphi\rightarrow\psi$, it follows that 
$A\models\psi$. That is, $\psi\in G(A)$. 

Suppose now that $\varphi$, $\psi\in G(A)$. This means that 
\begin{xitemize}
\xitem[buk-9-0] 
  $A\models\varphi$ and
  $A\models\psi$. 
\end{xitemize}
Hence we have $A\models\varphi\land\psi$. 
By \Claimof{buk-poP}, it follows that $(\varphi\land\psi)\in\poP$, that is,
$T\not\vdash\neg(\varphi\land\psi)$. Thus 
$\varphi$ and $\psi$ are compatible by \Claimof{buk-comp}.
\qedofClaim

\begin{Claim}
  \label{buk-generic}
  $G(A)$ is $\poP$-generic.
\end{Claim}
\prfofClaim
Working in $M$, suppose that $\Gamma$ is a maximal antichain in $\poP$. By 
\Claimof{buk-kappa-cc}, we have $\cardof{\Gamma}<\kappa$ and hence we have
$\llor\Gamma\in\calL_{\lambda}(\mu)$ and hence 
$\llor\Gamma\in\poP$: For $\varphi\in\Gamma$, since $\varphi\in\poP$ we have
$T\not\vdash\neg\varphi$ and $\vdash\varphi\rightarrow\llor\Gamma$. It follows
$T\not\vdash\llor\Gamma$. 

Moreover we have $T\vdash\llor\Gamma$: Otherwise
$\neg\llor\Gamma$ would be  
an element of $\poP$ incompatible with every $\varphi\in\Gamma$. A contradiction 
to the maximality of $\Gamma$. 

Hence $A\models\llor\Gamma$ and thus there is $\varphi\in\Gamma$ \st\
$A\models\varphi$. That is, $\varphi\in G(A)$. \qedofClaim\\
\qedofLemma\qedskip

The proof of \Thmof{bukovsky-thm} from \Lemmaof{buk-main-cl} relies on 
\Lemmaof{kappa-c-c} and the Axiom of 
Choice is involved both in the statement and the proof of 
\Lemmaof{kappa-c-c}. 

On the  
other hand, \Lemmaof{buk-main-cl} can be proved without assuming
the Axiom of Choice in $M$: It suffices to eliminate choice from the proof of
\Claimof{buk-generic}.\qedskip

\noindent
{\bf Proof of \bfClaimof{buk-generic} without the Axiom of Choice in $M$:} Working 
in $M$, suppose that $D$ is a dense subset of $\poP$. Then 
$A\models\llor D$: Otherwise we would have $T\not\vdash\llor D$. Since
\begin{xitemize}
\xitem[buk-10] 
  $T\vdash\llor D\leftrightarrow\llor g(D)$, 
\end{xitemize}
it follows that $T\not\vdash \llor g(D)$. Since
$\llor g(D)\in\calL_{\lambda}(\mu)$, this implies 
$\neg\llor g(D)\in\poP$. Since $D$ is dense in $\poP$ there is $\varphi_0\in D$ 
\st\ $T\vdash\varphi_0\rightarrow\neg\llor g(D)$. By \xitemof{buk-10}, it follows 
that $T\vdash\varphi_0\rightarrow\neg\llor D$. On the other hand, since
$\varphi_0\in D$ we have $T\vdash\varphi_0\rightarrow\llor D$. Hence we have
$T\vdash\neg\varphi_0$ which is a contradiction to $\varphi_0\in\poP$. 

Thus there is $\varphi_1\in D$ \st\ $A\models \varphi_1$, that is, $\varphi_1\in G(A)$. \\
\qedof{\Claimof{buk-generic} without \AC\ in $M$}\qedskip

The next corollary follows immediately from this remark:
\begin{Cor}
  \label{withoutAC}
  Work in \NBG. Suppose that $V$ is a model of\/ \ZFC\ and  $M$ is an inner model of $V$ 
  (of\/ \ZF) \st\ $M$ $\kappa$-globally covers $V$. If $V=M[A]$ for some set
  $A\subseteq\On$ then $V$ is a $\kappa$-c.c.\ set-generic extension of $M$. \qed
\end{Cor}

 We do not know if \Corabove\ is false
without the added assumption that $V$ is $M[A]$ for a set of
ordinals $A$.

More generally, it seems to be open if there is a characterisation of the 
set-generic extensions of an arbitrary model of ZF; or at least of such 
extensions given by partial orders which are well-ordered in the ground model. 

Grigorieff's theorem can be also obtained by a modification of the proof of 
\Thmof{bukovsky-thm}.  
\begin{Cor}[S.\ Grigorieff \cite{grigorieff}]
\label{grigorieff}
Suppose that $M$ is an inner model of a model $V$ of \ZFC\ and $V$ is a set-generic extension of $M$. Then any inner 
model $N$ of $V$ (of \ZFC) with $M\subseteq N$ is a set-generic extension of $M$ and hence 
definable in $V$.  Also, for such $N$, $V$ is a set-generic extension of $N$. 

If $V$ is $\kappa$-c.c.\ set-generic extension of $M$ in addition, then $N$ 
is a $\kappa$-c.c.\ set-generic extension of $M$ and $V$ is a $\kappa$-c.c.\ set-generic extension of $N$. \qed
\end{Cor}


Similarly to \Thmof{bukovsky-thm}, we can also characterize generic extensions obtained 
via a partial ordering of cardinality $\leq\kappa$. 

For $M$ and $V$ as above, we say that $V$ is {\it$\kappa$-decomposable} into $M$ if 
for any $a\in V$ with $a\subseteq M$, there are $a_i\in M$, $i\in\kappa$ \st\ 
$a=\bigcup_{i<\kappa}a_i$. 

\begin{Thm}
  Suppose that $V$ is a transitive model of \ZFC\ and $M$ an inner model of \ZFC\ 
  definable in $V$ and $\kappa$ is a cardinal in $M$. Then $V$ is a generic 
  extension of $M$ by a partial ordering in $M$ of size $\leq\kappa$ (in $M$) if 
  and only if $M$ $\kappa^+$-globally covers $V$ and $V$ is $\kappa$-decomposable 
  into $M$. 
\end{Thm}
\prf If $V$ is a generic extension of $M$ by a generic filter $G$ over a partial 
ordering $\poP\in M$ of size $\leq\kappa$ (in $M$) then $\poP$ has 
the $\kappa^+$-c.c.\ and hence $M$ $\kappa^+$-globally covers $V$ by 
\Thmof{laver-thm}. $V$ is $\kappa$-decomposable into $M$ since, for any $a\in V$ 
with $a=\dota^G$, we have
$a=\bigcup\setof{\setof{m\in M}{p\forces{\poP}{m\in\dota}}}{p\in G}$. 

Suppose now that $M$ $\kappa^+$-globally covers $V$ and $V$ 
is $\kappa$-decomposable into $M$. By \Thmof{bukovsky-thm}, there is 
a $\kappa^+$-c.c.\ 
partial ordering $\poP$ in $M$ and a $\poP$-generic filter $G$ over $M$ \st\
$V=M[G]$. \Wolog, we may assume that $\poP$ consists of the positive elements of a 
complete Boolean algebra $\BAB$ (in $M$). 

By $\kappa$-decomposability, $G$ can be decomposed into $\kappa$ sets $G_i\in M$,
$i<\kappa$. \Wolog, we may assume that $\bbbone_\poP$ forces this fact. So 
letting $\dotG$ be the standard name of $G$ and $\dotG_i$, $i<\kappa$ be names of 
$G_i$, $i<\kappa$ respectively, we may assume
\begin{xitemize}
\xitem[buk-9-1] 
  $\forces{\poP}{\dotG=\bigcup_{i<\kappa}\dotG_i}$.
\end{xitemize}
Working in $M$, let $X_i\subseteq\poP$ be a maximal pairwise incompatible set of 
conditions $p$ which decide $\dotG_i$ to be $G_{i,p}\in M$ for each $i<\kappa$. 
By the $\kappa^+$-c.c.\ of $\poP$, we have $\cardof{X_i}\leq\kappa$.
Clearly, we have $p\leq_\poP\prod^\BAB G_{i,p}$ for all $i<\kappa$ and $p\in X_i$.
Let $\poP'=\bigcup\setof{X_i}{i<\kappa}$. Then $\cardof{\poP'}\leq\kappa$. 
\begin{Claim}
  $\poP'$ is dense in $\poP$. 
\end{Claim}
\prfofClaim
Suppose $p\in\poP$. Then there is $q\leq p$ \st\ $q$ decides some $\dotG_i$ to be 
$G_{i,q}$ and $p\in G_{i,q}$. Let $r\in X_i$ be compatible with $q$. Then we have
$r\leq_\poP \prod^\BAB G_{i,r}=\prod^\BAB G_{i,q}\leq p$. 
\qedofClaim\qedskip\\
Thus $V$ is a $\poP'$-generic extension over $M$. \qedofThm\qedskip

\section{A Formal deductive system for $\calL_{\infty}(\mu)$ }
\label{ded-sys}
In the proof of \Lemmaof{buk-main-cl}, we used a formal deductive system of
$\calL_{\infty}(\mu)$ without specifying exactly which system we are using. 
It is enough to  
consider a system of deduction which contains all logical axioms we used in the course of the proof 
together with modus ponens and some infinitary deduction rules like: 
\[
\begin{array}{c}
  \varphi_i\rightarrow\psi,\ \ \ i\in I\\[-0.4em]
  \rule{28ex}{0.4pt}\\
  \llor\setof{\varphi_i}{i\in I}\rightarrow\psi
\end{array}
\]\noindent
What we need for such a system is that its correctness and upward absoluteness hold 
while we do not make use of any version of completeness of the system.

\finalmod{Formal deduction systems for infinitary logics have been studied 
  extensively in 1960s and 1970s, see e.g. \cite{karp}, \cite{keisler}, 
  \cite{takeuti}. }
Nevertheless, to be concrete, we shall introduce below such a deductive system \sfS\ for
$\calL_{\infty}(\mu)$. \medskip

\finalmod{One peculiar task for us here is that we have to make our deduction 
  system \sfS\ \st\ \sfS\ does not rely on \AC\ so that we can apply it in an 
  inner model $M$ which does not necessarily satisfy \AC\ to obtain \Corof{withoutAC}.}

Recall that we have introduced $\calL_{\infty}(\mu)$ as the smallest class 
containing the sets $\pairof{\alpha,0}$, $\alpha\in\mu$ as the codes of the 
prediactes ``$\alpha\in\dot{A}$'' for $\alpha\in\mu$ and closed \wrt\
$\pairof{\varphi,1}$ for $\varphi\in\calL_{\infty}(\mu)$ and $\pairof{\Phi,2}$ 
for all sets $\Phi\subseteq\calL_{\infty}(\mu)$ where $\pairof{\varphi,1}$ and
$\pairof{\Phi,2}$ 
represent $\neg\varphi$ and $\llor\Phi$ respectively. Here, to be more precise about the 
role of the infinite conjunction we add the infinitary logical 
connective $\lland$, and assume that $\lland\Phi$ is coded by $\pairof{\Phi,3}$ and 
thus $\calL_{\infty}(\mu)$ is also closed \wrt\ 
$\pairof{\Phi,3}$ for all sets $\Phi\subseteq\calL_{\infty}(\mu)$. 

The axioms of \sfS\ consist of the following formulas:
\begin{description}
\item[(A1)] $\varphi(\varphi_0,\varphi_1\ctentenc\varphi_{n-1})$\medskip\\
  for each 
  tautology $\varphi(A_0,A_1\ctentenc A_{n-1})$ of (finitary) propositional logic 
  and $\varphi_0$, $\varphi_1$\ctentenc $\varphi_{n-1}\in\calL_{\infty}(\mu)$; 
\item[(A2)] $\varphi\rightarrow\llor\Phi$\ \ and\ \ 
  $\lland\Phi\rightarrow\varphi$\medskip\\
  for any set $\Phi\subset\calL_{\infty}(\mu)$ and $\varphi\in\Phi$; 
\item[(A3)] $\neg(\lland\Phi)\leftrightarrow
  \llor\setof{\neg\varphi}{\varphi\in\Phi}$ and\\ 
  $\neg(\llor\Phi)\leftrightarrow
  \lland\setof{\neg\varphi}{\varphi\in\Phi}$\medskip\\
  for any set 
  $\Phi\subseteq\calL_{\infty}(\mu)$; and 
\item[(A4)] $\varphi\land(\llor\Psi)\leftrightarrow
  \llor\setof{\varphi\land\psi}{\psi\in\Psi}$ and\\ 
  $\varphi\lor(\lland\Psi)\leftrightarrow
  \lland\setof{\varphi\lor\psi}{\psi\in\Psi}$\medskip\\
  for any $\varphi\in\calL_{\infty}(\mu)$ and any set $\Psi\subseteq\calL_{\infty}(\mu)$. 
\end{description}

Deduction Rules: 

\begin{description}
  \item[(Modus Ponens)]\ \ \   
    $
    \begin{array}{c}
      \ssetof{\varphi, \varphi\rightarrow\psi}\\[-0.4em]
      \rule{18ex}{0.4pt}\\
      \psi
    \end{array}
    $
  \item[(R1)]\ \ \  
    $
    \begin{array}{c}
      \setof{\varphi\rightarrow\psi}{\varphi\in\Phi}\\[-0.4em]
      \rule{28ex}{0.4pt}\\
      \llor\Phi\rightarrow\psi
    \end{array}
    $\quad {\bf (R2)}\ \ \ 
    $
    \begin{array}{c}
      \setof{\varphi\rightarrow\psi}{\psi\in\Psi}\\[-0.4em]
      \rule{28ex}{0.4pt}\\
      \varphi\rightarrow\lland\Psi
    \end{array}
    $
\end{description}

A proof of $\varphi\in\calL_{\infty}(\mu)$ from
$\Gamma\subseteq\calL_{\infty}(\mu)$ is a labeled tree $\pairof{\sfT,f}$ \st\ 
\begin{xitemize}
\xitem[prf-0] $\sfT=\pairof{\sfT,\leq}$ is a tree growing upwards with its root $r_0$ 
  and \finalmod{$\sfT$ with $(\leq)^{-1}$ is well-founded};
\xitem[prf-1] $\mapping{f}{\sfT}{\calL_{\infty}(\mu)}$;
\xitem[prf-1-0] $f(r_0)=\varphi$;
\xitem[prf-2] if $t\in\sfT$ is a maximal element then either $f(t)\in\Gamma$ or $t$ 
  is one of the axioms of $\sfS$;
\xitem[prf-3] if $t\in\sfT$ and $P\subseteq\sfT$ is the set of all immediate 
  successors of $t$, then 
  \[ \begin{array}{c}
      \setof{f(p)}{p\in P}\\[-0.4em]
      \rule{20ex}{0.4pt}\\
      f(t)
    \end{array}
  \]\noindent
  is one of the deduction rules.
\end{xitemize}

\finalmod{%
We have to stress here that, in \xitemof{prf-3}, we do not assume that the 
function $f$ is one-to-one since otherwise we have to choose a proof for each 
formula in the set in the premises of \itemof{R1} and \itemof{R2}. 
Thus, for example, we can deduce $T\vdash\lland\Phi$ in \sfS\ from 
$T\vdash\varphi$ for all $\varphi\in\Phi$ without appealing to \AC. }

Now the proof of the following is an easy exercise:
\begin{Prop}
  \assertof{1} For any $B\subseteq\mu$, $T\subseteq\calL_{\infty}(\mu)$ and
  $\varphi\in\calL_{\infty}(\mu)$, if $T\vdash\varphi$ and $B\models T$, then 
  we have $B\models\varphi$.\smallskip

  \assertof{2} For transitive models $M$, $N$ of\/ \ZF\ \st\ $M$ is an inner model 
  of $N$, if 
  $M\models\mbox{``\/}\pairof{\sfT,f}\mbox{ is a proof of }\varphi\mbox{ in }
  \calL_{\infty}(\mu)\mbox{''}$, then 
  \begin{xitemize}
  \item[] 
    $N\models\mbox{``\/}\pairof{\sfT,f}\mbox{ is a proof of }\varphi\mbox{ in }
    \calL_{\infty}(\mu)\mbox{''}$. 
  \end{xitemize}
\end{Prop}
\prf \assertof{1}: By induction on cofinal subtrees of a fixed proof
$\pairof{\sfT,f}$ of $\varphi$.  
\assertof{2}: Clear by definition. \qedofProp\qedskip

An alternative setting to the argument by means of a deductive system is to make use of 
the following definition of  
$M\models\mbox{``}\Gamma\vdash \varphi\mbox{''}$ in the proof of \Lemmaof{buk-main-cl}: 
\begin{xitemize}
\item[] 
  $M\models\mbox{``}\Gamma\vdash\varphi\mbox{''}$ iff
  for any $B\subseteq\mu$ in some set-forcing extension $M[G]$ of $M$,
  $M[G]\models B\models\psi$ for all $\psi\in\Gamma$ always  
  implies $M[G]\models B\models\varphi$. 
\end{xitemize}
Note that this is definable in $M$ using the forcing relation definable on $M$. It 
remains to verify that  
this notion has the desired degree of absoluteness. Actually we can easily prove 
the full absoluteness, that is, if $N$ is a transitive 
model containing $M$ 
with the same ordinals as those of $M$ then, for $\Gamma$, $\varphi\in M$ with
$M\models\Gamma\subseteq\calL_{\infty}(\mu)$ and
$M\models\varphi\in\calL_{\infty}(\mu)$, 
$\Gamma\vdash\varphi$ holds in $M$ iff
$\Gamma\vdash\varphi$ holds in $N$.  

First suppose that $B\subseteq\mu$ is 
a set of ordinals in a set-generic 
extension $N[G]$ of $N$ \st\ $B$ witnesses the failure of $\Gamma\vdash\varphi$ 
in $N$. Let $x$  
be a real which is 
generic over $N$ for the L\'evy collapse of a sufficiently large $\nu$ to $\omega$ 
\st\ $\Gamma$ and $\mu$ become countable in the generic extension $N[x]$. 
Then $x$ is also L\'evy generic over $M$ and $M[x]$ is a submodel of $N[x]$. By
L\'evy Absoluteness, it follows that that there exists $B'\subseteq\mu$ in
$M[x]$ which also 
witnesses the failure of 
$\Gamma\vdash\varphi$ in $M$. 

Conversely, suppose that $\Gamma\vdash\varphi$ 
holds in $N$ and let 
$B\subseteq\mu$ be a set of ordinals in a set-generic extension $M[G]$ of $M$ 
\st\ $B$ witnesses the failure of 
$\Gamma\vdash\varphi$ in $M$. Then $B$ also belongs to an
extension of $M$ which is generic for the L\'evy collapse of sufficently large
$\nu$ to $\omega$;  
choose a condition  
$p$ in this forcing which forces the existence of such a $B$. Now if $x$ is 
L\'evy-generic over 
$N$ and contains the condition $p$, we see that there is a counterexample to
$\Gamma\vdash\varphi$  
in $N$ witnessed in $N[x]$, contrary to our assumption.

With both of the interpretations of $\vdash$ we can check that the arguments 
in \sectionof{laver-bukovsky} go through. 

\section{An axiomatic framework for the set-generic multiverse}
\label{multi-sys}
In this section, \finalmod{we consider some possible axiomatic treatments of 
  the set-generic multiverse. Such axiomatic treatments are also
  discussed e.g.\ in \cite{gitman-hamkins}, \cite{steel},
  \cite{vaananen}.} 
We introduce a conservative extension \MZFC\ of \ZFC\ in which 
we can treat the multiverse of set-generic extensions of models of \ZFC\ as a 
collection of countable transitive models. \finalmod{This system or some 
further extension of it (which can possibly  
also treat tame class forcings) may be used as a basis 
for direct formulation of statements concerning the multiverse. }

\newcommand{\ttv}{{\mathtt v}}
The language $\calL_\MZF$ of the axiom system \MZFC\ consists of 
the $\epsilon$-relation symbol `$\in$', and a constant symbol `$\ttv$' which should 
represent the countable transitive ``ground model''.

The axiom system \MZFC\ consists of
\begin{xitemize}
\xitem[mzfc-0] all axioms of \ZFC;
\xitem[mzfc-1]
  ``$\ttv\mbox{ is a countable transitive set}$'';
\xitem[mzfc-2] ``$\ttv\models\varphi$'' for all axioms $\varphi$ of \ZFC;
\end{xitemize}

By \xitemof{mzfc-0}, \MZFC\ proves 
the (unique) existence of the closure $\calM$ of ``$\ssetof{\ttv}$'' under 
forcing extension and definable ``inner model'' of ``\ZF'' (here `\ZF'\ is set in 
quotation marks since we can only argue in metamathematics that such ``inner model'' 
satisfies each instance of replacement). Note that $\calM\subseteq\calH_{\aleph_1}$. 
Here ``inner model'' is actually phrased in $\calL_\ZF$ as ``transitive almost universal 
subset closed under G\"odel operations''. If we had $\ttv\models\ZFC$, we 
would have $w\models\ZF$ for any inner model $w$ of $\ttv$ in this sense by Theorem 13.9 in  
\cite{millennium-book}.  In \MZFC, however, we have only $\ttv\models\varphi$ for 
each axiom $\varphi$ of \ZFC\ (in the meta-mathematics). Nevertheless, for 
all such ``inner model'' $w$ and hence for all $w\in\calM$, 
we have $w\models\varphi$ for all 
axiom $\varphi$ of \ZF\ 
by the proof 
of  Theorem 13.9 in \cite{millennium-book} and the Forcing Theorem.  Apparently, 
this is enough to consider $\calM$ in this framework as the set-generic multiverse. 

Similarly, we can also start from any extension of \ZFC\ (e.g.\, with 
some additional large cardinal axiom) and make 
$\calM$ closed under some more operations such as some well distinguished class of 
class forcing extensions. 

The following theorem shows that we do not increase the 
consistency strength by moving from \ZFC\ to \MZFC. 

\begin{Thm}
  \MZFC\ is a conservative extension of \ZFC: for any sentence $\psi$ in
  $\calL_\ZF$, we have $\ZFC\vdash\psi$\ \ $\Leftrightarrow$\ \ $\MZFC\vdash\psi$.
  In particular, \MZFC\ is equiconsistent with \ZFC.
\end{Thm}
\prf ``$\Rightarrow$'' is trivial. 

For ``$\Leftarrow$'', suppose that $\MZFC\vdash\psi$ for a formula $\psi$ in
$\calL_\ZF$. Let $\calP$ be a proof of $\psi$ from \MZFC\ and let $T$ be the 
finite fragment of \ZFC\ consisting of all axioms $\varphi$ of \ZFC\ \st\
$\ttv\models\varphi$ appears in $\calP$. Let $\Phi(x)$ be the formula in 
$\calL_\ZF$ saying
\begin{xitemize}
\item[] ``$x$ is a countable transitive set and $x\models\lland T$''.
\end{xitemize}
By the Deduction Theorem, we can recast $\calP$ to a proof of
$\ZFC\vdash\forall x(\Phi(x)\rightarrow\psi)$. On the other hand we have
$\ZFC\vdash\exists x\Phi(x)$ (by the Reflection Principle, Downward 
L\"owenheim-Skolem Theorem and Mostowski's Collapsing Theorem). Hence we obtain a proof of
$\psi$ from \ZFC\ alone. 
\qedofThm
\qedskip

It may be a little bit disappointing if each set-theoretic universe in the multiverse seen
from the ``meta-universe'' is merely a countable set.
Of course 
if $M$ is an inner model of a model $W$ of \ZFC\ (i.e. $M$ is a model which is a 
transitive class $\subseteq W$ and $M$, $W\models\ZFC$) there are always \po\
$\poP$ in $M$ for which  
there is no $(M,\poP)$-generic set in $W$ (e.g. any \po\ collapsing a cardinal of 
$W$ cannot have its generic set in $W$). 

However, if we are content with a meta-universe which is not a model of full 
\ZFC, we can work with the following setting where each of the ``elements'' of 
the set-generic  
multiverse is an inner model of a meta-universe: starting from a model $V$ of \ZFC\ 
with an inaccessible cardinal $\kappa$, we generically extend it to $W=V[G]$ by 
L\'evy collapsing $\kappa$ to $\omega_1$. Letting $M=\calH(\kappa)^V$, we have
$M\models\ZFC$ and $M$ is an inner model of
$W=\calH(\kappa)^{V[G]}=\calH(\omega_1)^{V[G]}$. $W\models\ZFC-{}$ the Power Set 
Axiom and for any \po\ $\poP$ in $M$ there is a $(M,\poP)$-generic set in $W$. 
Thus an NBG-type theory of $W$ with a new unary predicate corresponding to $M$ can 
be used as a framework of the theory for the set-generic multiverse (which is 
obtained by considering all the set-generic grounds of $M$, and then all the 
set generic extensions of them, etc.) as a ``class'' of classes in $W$. 
\finalmod{A setting similar to this idea was also discussed in \cite{steel}.}

\section{Independent buttons}
\label{ind-buttons}
The multiverse view sometimes highlights problems which would be never 
asked in the conventional context of forcing constructions. The existence of 
infinitely many independent buttons which arose in connection with the characterization of the
modal logic of the set-generic multiverse (see \cite{hamkins-loewe}) is one such question. 

A sentence $\varphi$ in $\calL_\ZF$ is said to be a {\it button\/} (for set-genericity) if 
any set-generic extension $V[G]$ of the ground model $V$ has a further set-generic
extension $V[G][H]$ such that $\varphi$ holds in all set-generic extensions of $V[G][H]$.
Let us say that a button $\varphi$ is pushed in a set-generic extension  
$V[G]$ if $\varphi$ holds in all further set-generic extensions $V[G][H]$ of $V[G]$ (including 
$V[G]$ itself). 

Formulas $\varphi_n$, $n\in\omega$ are {\it independent buttons}, if $\varphi_n$, 
$n\in\omega$ are unpushed buttons and for any set-generic extension $V[G]$ of the ground model
$V$ and any $X\subseteq\omega$ in $V[G]$, 
\begin{xitemize}
\xitem[button-0] 
  if
  $\setof{n\in\omega}{V[G]\models\varphi_n\mbox{ is pushed}}\subseteq X$ then 
  there is a set-generic  
  extension $V[G][H]$ \st\ $\setof{n\in\omega}{V[G][H]\models\varphi_n\mbox{ is pushed}}=X$. 
\end{xitemize}

In \cite{hamkins-loewe}, it is claimed that formulas 
$b_n$, $n\in\omega$ form an infinite set of independent buttons over $V=L$ 
where $b_n$ is a formula asserting: ``${\omega_n}^L$ is not a cardinal''. 
This is used to prove that the principles of forcing 
expressible in the modal logic of 
the set-theoretic multiverse as a Kripke frame where modal operator $\Box$ 
is interpreted as: 
\begin{xitemize}
\xitem[button-1] 
  $M\models\Box\varphi$ $\Leftrightarrow$ in all 
  set-generic extensions $M[G]$ of $M$ we have $M[G]\models\varphi$  
\end{xitemize}
coincides with the modal theory {\sf S4.2} (Main Theorem 6 in \cite{hamkins-loewe}). 

Unfortunately, it seems that there is no guarantee that \xitemof{button-0} holds  
in an arbitrary set-generic extension $V[G]$ for these $b_n$, $n\in\omega$. 

In the following, we introduce an alternative set of infinitely many formulas 
which are actually independent buttons for any ground model of \ZFC $+$ ``\GCH\ 
below $\aleph_\omega$'' $+$ ``$\aleph_n=\aleph_n^L$ for all $n\in\omega$''   
which can be used as $b_n$, $n\in\omega$ in \cite{hamkins-loewe}. 

We first note that, for Main Theorem 6 in \cite{hamkins-loewe} we actually need 
only the existence of an arbitrary finite number of independent buttons. In the case 
of $V=L$ the following formulas can be used for this: Let $\psi_n$ 
be the statement that $\aleph_n^L$ is a cardinal and the $L$-least $\aleph_n^L$-Suslin tree 
$T^L_n$ in $L$ (i.e., the $L$-least normal tree of height 
$\aleph_n^L$ with no antichain of size $\aleph_n^L$ in $L$) is 
still $\aleph_n^L$-Suslin. If $M$ is a  
set-generic (or arbitrary)
extension of $L$ in which the button $\neg\psi_n$ has not been pushed, then by forcing with 
$T_n^L$ over $M$ we push this button and do not affect any of the other unpushed 
buttons $\neg\psi_m$, $m\neq n$, 
as this forcing is $\aleph_n$-distributive and has size
$\aleph_n$. Rittberg \cite{rittberg} also found independent buttons under $V=L$. 

Now we turn to a construction of infinitely many independent buttons for which we 
even do not need the existence of Suslin trees.
For $n\in\omega$, let $\varphi_n$ be the statement:
\begin{xitemize}
\xitem[button-2] there is an injection from ${\aleph_{n+2}}^L$ to
  $\psof{{\aleph_n}^L}$. 
\end{xitemize}
Note that $\varphi_n$ is pushed in a set-generic extension $V[G]$ if and only if it 
holds in $V[G]$. 
Thus $\varphi_n$ for each $n\in\omega$ is a button provided that $\varphi_n$ 
does not hold in the ground model. 
We show that these $\varphi_n$, $n\in\omega$ are independent buttons (over any 
ground model where they are unpushed --- e.g., when $V=L$). 

Suppose that we are working in some model $W$ of \ZFC. In $W$, let
$A=\setof{n\in\omega}{\Box\varphi_n\mbox{ holds}}$ and $B\subseteq\omega$ be 
arbitrary with $A\subseteq B$. It is enough to prove the following

\begin{Prop}
  We can force (over $W$) that $\varphi_n$ holds for all $n\in B$ and 
$\neg\varphi_n$ for all $n\in\omega\setminus B$. 
\end{Prop}
\prf
In $W$, let $\kappa_n=\cardof{{\aleph_n}^L}$ for $n\in\omega$. We use the 
notation of \cite{kunen} on the partial orderings with partial functions and 
denote with $\Fn(\kappa,\lambda,\mu)$ the set of all partial functions from
$\kappa$ to $\lambda$ with cardinality $<\mu$ ordered by reverse inclusion. By
$\Delta$-System Lemma, it is easy to see that $\Fn(\kappa,\lambda,\mu)$ has the
$(\lambda^{<\mu})^+$-c.c. 
Let 
\begin{xitemize}
\xitem[button-3] $\poP_n=\left\{\,
  \begin{array}{@{}ll}
    \Fn(\kappa_{n+2}, 2,\kappa_n) &\mbox{if }n\in B\setminus A\\
    \bbbone &\mbox{otherwise.}
  \end{array}
  \right.$
\end{xitemize}
Let $\poP=\prod_{n\in\omega}\poP_n$ be the full support product of $\poP_n$,
$n\in\omega$. Then we clearly have $\forces{\poP}{\varphi_n}$ for all $n\in B$. 
Thus to show that $\poP$ creates a generic extension as desired, it is enough 
to show that $\forces{\poP}{\neg\varphi_n}$ for all $n\in\omega\setminus B$. 

Suppose that 
\begin{xitemize}
\xitem[button-3-0] 
  $n\in\omega\setminus B$. 
\end{xitemize}
Then we have 
\begin{xitemize}
\xitem[button-4] $\poP_n=\bbbone$. 
\end{xitemize}
Since $\varphi_n$ does not hold in $W$, we have
$\kappa_n<\kappa_{n+1}<\kappa_{n+2}$ and $2^{\kappa_n}=\kappa_{n+1}$ in $W$. By 
\xitemof{button-4}, $\poP$ factors as 
$\poP\sim\poP(<n)\times\poP(>n)$ where $\poP(<n)=\prod_{k<n}\poP_k$ and
$\poP(>n)=\prod_{k>n}\poP_k$. 

We show that both $\poP(>n)$ and $\poP(<n)$ over $\poP(>n)$ do not add any 
injection from $\kappa_{n+2}$ into $\psof{\kappa_n}$. 

$\poP(>n)$ is $\kappa_{n+1}$-closed. Thus it does not add any new subsets of
$\kappa_n$. So if it added an injection from $\kappa_{n+2}$ 
into $\psof{\kappa_n}$ then it would collapse the cardinal $\kappa_{n+2}$. Since 
$\poP(>n)$ further factors as $\poP(>n)\sim\poP_{n+1}\times\prod_{k>n+1}\poP_k$ 
and $\prod_{k>n+1}\poP_k$ is $\kappa_{n+2}$-closed the only way $\poP(>n)$ could 
collapse $\kappa_{n+2}$ would be if $\poP_{n+1}$ did so. But then, since 
$\poP_{n+1}$ has the $(2^{<\kappa_{n+1}})^+$-c.c.\ with
$(2^{<\kappa_{n+1}})^+=(2^{\kappa_n})^+$, we would have
$2^{\kappa_n}\geq\kappa_{n+2}$. This is a contradiction to the choice 
\xitemof{button-3-0} of $n$. So $\poP(>n)$ forces $\varphi_n$ to fail.

In the rest of the proof, we work in $W^{\poP(>n)}$ and show that $\poP(<n)$ does 
not add any injection from $\kappa_{n+2}$ into $\psof{\kappa_n}$. Note that, by
$\kappa_{n+1}$-closedness of $\poP(>n)$, we have
$\Fn(\kappa_{m+2},2,\kappa_m)^W=\Fn(\kappa_{m+2},2,\kappa_m)^{W^{\poP(>n)}}$ for
$m<n$. 

We have the following two cases:

\noindent
{\bf Case I. } $n-1\in A\cup(\omega\setminus B)$. Then $\poP(<n)\sim\poP(<m)$ for 
some $m<n$ and $\poP(<m)$ has the $(2^{\kappa_{m-1}})^+$-c.c.\ with
$(2^{\kappa_{m-1}})^+\leq\kappa_n$. \smallskip

\noindent
{\bf Case I{}I. } $n-1\in B\setminus A$. Then $2^{<\kappa_{n-1}}=\kappa_n$ and 
$\poP(<n)$ has the $\kappa_{n+1}$-c.c.

In both cases the partial ordering $\poP(<n)$ has $\kappa_{n+1}$-c.c.\ and hence 
the cardinals $\kappa_{n+1}$ and $\kappa_{n+2}$ are preserved. Since $\poP(<n)$ 
has at most cardinality $2^{\kappa_{n-1}}\cdot\kappa_{n+1}=\kappa_{n+1}$, it 
adds at most ${\kappa_{n+1}}^{\kappa_n}=\kappa_{n+1}$ new subsets of $\kappa_n$ 
and thus the size of $\psof{\kappa_n}$ remains unchanged. This shows that
$\forces{\poP}{\neg\varphi_n}$. \\
\qedofProp

\end{document}